\documentclass[MSNbibl,number,citesort,seceqn,dvips]{arxbj}
\usepackage{mathrsfs}

\aid{0}
\volume{20}
\issue{2}
\pubyear{2014}
\firstpage{990}
\lastpage{1005}
\doi{10.3150/13-BEJ513} %kopijuoti is PTS

\makeatletter
\newcommand{\mrmd}{\,\mathrm{d}}

\newcommand{\rrVert}{\Vert}
\newcommand{\rrvert}{\vert}
\newcommand{\llVert}{\Vert}
\newcommand{\llvert}{\vert}

\newcommand{\xrightarrow}[1]{\stackrel{#1}{\longrightarrow}}

\newtheorem{theorem}{Theorem}[section]

\newremark{remark}{Remark}[section]
\newremark{example}{Example}[section]

\makeatother

\begin{document}
\begin{frontmatter}

\title{Conditions for convergence of random coefficient AR(1)
processes and perpetuities in higher dimensions}
\runtitle{Convergence of random coefficient AR(1) processes}

\begin{aug}
%%%% inicialai - be tarpu
\author{\fnms{Torkel} \snm{Erhardsson}\corref{}\ead[label=e1]{torkel.erhardsson@liu.se}}% \and
\runauthor{T. Erhardsson} %% auto
\address{Department of Mathematics, Link\"oping University, S-581 83
Link\"oping,
Sweden.\\ \printead{e1}}
\end{aug}

% HISTORY:
\received{\smonth{7} \syear{2012}}
\revised{\smonth{1} \syear{2013}}

% ABSTRACT
%
\begin{abstract}
A $d$-dimensional RCA(1) process is a generalization of the
$d$-dimensional AR(1) process, such that the coefficients
$\{M_t;t=1,2,\ldots\}$ are i.i.d. random matrices. In the case $d=1$,
under a nondegeneracy condition, Goldie and Maller gave necessary and
sufficient conditions for the convergence in distribution of an RCA(1)
process, and for the almost sure convergence of a closely related sum
of random variables called a \emph{perpetuity}. We here prove that
under the condition $\llVert {\prod_{t=1}^nM_t}\rrVert
\xrightarrow{\mathrm{a.s.}}0$ as
$n\to\infty$, most of the results of Goldie and Maller can be extended
to the case $d>1$. If this condition does not hold, some of their
results cannot be extended.
\end{abstract}

% KEYWORDS
% visi is mazosios raides ir pagal abecele
%

\end{frontmatter}

%s1 #&#
\section{Introduction} \label{Sintroduction}
In this paper, we consider a discrete time stochastic process called
the \emph{$d$-dimensional $\operatorname{RCA}(1)$ process}, or \emph{random coefficient
autoregressive process of order 1}, which is a generalization of the
$d$-dimensional AR(1) process. We also consider a closely related
infinite sum of $d$-dimensional random variables, called a
\emph{perpetuity}. Since the appearance of \cite {Kes73}, different
aspects of the RCA(1) process and the perpetuity have been studied by
many authors; see, for example,
\mbox{\cite{Ver79,Gri80,Gri81,Bra86,BP92,EG94,GM00,DGL04,AI09}} and the
references therein. In the present work, we will focus on conditions
for convergence in distribution of the RCA(1) process, and for almost
sure convergence of the perpetuity.

For each positive integer $p$, the $d$-dimensional RCA($p$) process is
defined as follows. Let $\{(M_{t,1},\ldots,M_{t,p});t=1,2,\ldots\}$ be
an i.i.d. sequence of $p$-tuples of random matrices of dimension
$d\times
d$ (the \emph{coefficients}); let
$\{Z_t;t=1,2,\ldots\}$ be i.i.d. $d$-dimensional random variables
independent of the random matrices (the \emph{error} variables); and
let $Z_0$ be a $d$-dimensional random variable
independent of everything else (the \emph{initial state}). Define the
$d$-dimensional RCA($p$) process
$\{X_t;t=1,2,\ldots\}$ by
\[
X_0=Z_0;\qquad X_t = \sum
_{i=1}^{p\wedge t}M_{t,i}X_{t-i} +
Z_{t}\qquad\forall t=1,2,\ldots.
\]
If the distribution of $(M_{1,1},\ldots,M_{1,p})$ is degenerate at a
constant matrix $p$-tuple, the usual $d$-dimensional AR($p$) process is
obtained. However, for
the AR($p$) process it is often assumed that the error
variables have finite second moments. Here, we make no
such assumption.

The AR($p$) process was originally proposed as a statistical model for time
series, and it is today one of the most widely used such models. The
RCA($p$) process was first considered as a statistical model in
\cite{And76}. A much studied problem is under what conditions on
the coefficients there exists an RCA($p$) or AR($p$) process which is
(wide sense) stationary. For some answers to this problem,
and more information on these processes, see \cite{BD91,BP92,BL10,And76,NQ82}, and the references
therein.

The case $p=1$ has received special attention, since the RCA($1$)
process is easily
seen to be a \emph{Markov chain} on the state space
$(\mathbb{R}^d,\mathscr{R}^d)$. For such a process, it is natural to
ask under what conditions on the error
variables and the random coefficient the process is (Harris)
recurrent, positive, or convergent in distribution. For some partial
answers to these questions, see \cite{MT09} and the references
therein. See also
\cite{FT89} for a connection between RCA(1)
processes and Dirichlet processes; this connection was exploited in
\cite{Erh08} to construct a new method to carry out Bayesian inference
for an unknown finite measure, when a number of integrals with respect
to this measure has been observed.

The \emph{perpetuity} associated with a $d$-dimensional
RCA(1) process is defined as the almost sure limit (if the limit
exists) of the $d$-dimensional random sequence
$\{V_t;t=1,2,\ldots\}$, defined by:
\[
V_t = \sum_{i=1}^{t}\prod
_{j=1}^{i-1}M_jZ_{i}
\qquad\forall t=1,2,\ldots.
\]
The existence of the perpetuity is closely related to the convergence
in distribution of the $d$-dimensional RCA(1) process. In particular,
it is shown in Section \ref{Smaintheorem} that if $\llVert
{\prod_{t=1}^nM_t}\rrVert \xrightarrow{\mathrm{a.s.}}0$ as $n\to\infty$
(a condition to be called C0 below), then the two convergence
statements are equivalent. Moreover, in the case $d=1$, if
$\mathbb{P}(Z_1=0)<1$, it was shown in \cite{GM00} that the existence
of the perpetuity implies C0.

The main result in \cite{GM00}, their Theorem 2.1, is a \emph{complete
solution} in the case $d=1$ to the problem: under what conditions on
the error variables and the random coefficients does the perpetuity
exist? Five different conditions on the random
variables are given, which, if
$\mathbb{P}(Z_1=0)<1$, are shown to be equivalent, and to imply both
the existence of the perpetuity, and C0. Furthermore, it is shown that
under a certain
``nondegeneracy'' condition, the five conditions are necessary
for the convergence in distribution of the associated RCA(1) process.

The main result of the present paper,
Theorem \ref{Tmultidimperpetuities}, is a generalization of most of
Theorem 2.1
in \cite{GM00} to the case $d>1$. All except one of the conditions in the
latter theorem are considered. (It is unclear how the remaining
condition, which involves the
finiteness of a particular integral, should be generalized to the
case $d>1$, if indeed this is possible at all.)
It is shown that if C0 is assumed, the remaining conditions of Theorem
2.1 are equivalent, and
imply the existence of the perpetuity. However, contrary to the
case $d=1$, the conditions do not imply C0, and if C0 is not
assumed, they are not all equivalent. Similarly,
under C0, the existence of the perpetuity is equivalent to
the convergence in distribution of the associated $d$-dimensional
RCA(1) process; not so without C0.

The remaining part of the paper is structured as follows: in
Section \ref{Smaintheorem}, the main result is stated and
proven; in Section \ref{Scounterexamples}, some counterexamples and
special cases are collected; and Section \ref{Sfurtherwork}
contains some suggestions for future research.

%s2 #&#
\section{Main result and proof} \label{Smaintheorem}
Let $d$ be a positive integer. Denote by $\llvert \cdot\rrvert $ the
Euclidean norm on the space $\mathbb{R}^d$. Let $\mathbb{R}^{d\times
d}$ be the space of $d\times d$-matrices with elements in
$\mathbb{R}$, and denote by \mbox{$\llVert \cdot\rrVert $} the
matrix norm
induced by \mbox{$\llvert \cdot\rrvert $}, that is, $\llVert
A\rrVert  =
\max_{\llvert  x\rrvert =1}\llvert  Ax\rrvert $. (This
is known as the
\emph{spectral norm}, and is equal to the largest singular value of
$A$.) Denote by $I_d$ the identity
$d\times d$-matrix. The following notation will
be used for matrix products:
\[
\prod_{j=m}^nM_j = \cases{
M_mM_{m+1}\cdots M_n, &\quad if $m\leq n$;
\cr
I_d, &\quad if $m>n$.}
\]
In particular, $\prod_{j=m}^{n-1}M_{n-j} =
M_{n-m}M_{n-m-1}\cdots M_1$ for each $m<n$, and
$\prod_{j=m}^{n-1}M_{n-j} = I_d$ for each $m\geq n$. Lastly, by
convention a minimum over an empty set is defined as $\infty$.

%th2.1 #&#
\begin{theorem} \label{Tmultidimperpetuities}
Let $\{(M_t,Z_t);t=1,2,\ldots\}$ be i.i.d. random elements in
$(\mathbb{R}^{d\times d}\times\mathbb{R}^d,\mathscr{R}^{d\times
d}\times\mathscr{R}^d)$, and let $Z_0$ be a random element in
$(\mathbb{R}^d,\mathscr{R}^d)$
independent of $\{(M_t,Z_t);t=1,2,\ldots\}$. Define the random sequence
$\{X_t;t=1,2,\ldots\}$ by
\[
X_0=Z_0;\qquad X_t = M_tX_{t-1}
+ Z_{t}\qquad\forall t=1,2,\ldots.
\]
Under the condition \emph{C0}: \emph{$\llVert \prod_{t=1}^nM_t\rrVert \xrightarrow{\mathrm{a.s.}}0$} as
$n\to\infty$, the following are equivalent:
\begin{eqnarray*}
&&\mbox{\textup{\hphantom{ii}(i)}}\quad X_t \mbox{ converges in distribution}\qquad\mbox{as $t\to
\infty$};
\\
&&\mbox{\textup{\hphantom{i}(ii)}}\quad \sum_{t=1}^\infty\Biggl|\prod
_{j=1}^{t-1}M_jZ_{t}\Biggr|<
\infty\qquad\mbox{a.s.};
\\
&&\mbox{\textup{(iii)}}\quad \sum_{i=1}^t\prod
_{j=1}^{i-1}M_jZ_{i}
\mbox{ converges a.s.}\qquad \mbox{as $t\to\infty$};
\\
&&\mbox{\textup{\hspace*{0.7pt}(iv)}}\quad \prod_{j=1}^{t-1}M_jZ_{t}
\xrightarrow{\mathit{a.s.}}0 \qquad\mbox{as $t\to\infty$};
\\
&&\mbox{\textup{\hspace*{0.5pt}\hphantom{i}(v)}}\quad \sup_{t=1,2,\ldots}\Biggl|\prod
_{j=1}^{t-1}M_jZ_{t}\Biggr|<\infty\qquad
\mbox{a.s.};
\\
&&\mbox{\textup{\hspace*{0.5pt}(vi)}}\quad \sum_{t=1}^\infty
\mathbb{P}\Biggl(\min_{k=1,\ldots,t-1}\Biggl|\prod_{j=k}^{t-1}M_jZ_t\Biggr|>
x\Biggr)<\infty \qquad\forall x>0.
\end{eqnarray*}
\end{theorem}

%re2.1 #&#
\begin{remark} \label{Rimplicationswithoutnormconv}
Clearly, the implications
(ii)${}\Rightarrow{}$(iii)${}\Rightarrow{}$(iv)${}\Rightarrow{}$(v) remain valid
even if C0 does not hold, and, as will be seen from the proof, so does the
implication (iv)${}\Rightarrow{}$(vi). It will
be shown in Example \ref{Rvtovicounterexample} that the
implication (v)${}\Rightarrow{}$(vi) need not hold if C0 does not hold. On
the other hand, in the case $d=1$, it was shown in
\cite{GM00} that if $\mathbb{P}(Z_1=0)<1$, then
(vi) implies C0, and if
also $\mathbb{P}(|M_1|=1)<1$, then (v) implies C0; see
Example \ref{RC0necessarycounterexample} below. -- The almost sure
limit of the sum in (iii) is called a
\emph{perpetuity}. Hence, (iii) is the statement that the perpetuity
exists.
\end{remark}

\begin{pf*}{Proof of Theorem \ref{Tmultidimperpetuities}}
(iii)${}\Rightarrow{}$(i). As is easily shown by induction, we can write
\[
X_t = \sum_{i=0}^{t-1}\prod
_{j=0}^{i-1}M_{t-j}Z_{t-i} +
\prod_{j=0}^{t-1}M_{t-j}Z_0
\qquad\forall t=1,2,\ldots.
\]
Replacing $(M_{t-i},Z_{t-i})$ by $(M_{i+1},Z_{i+1})$ for $i=0,1,\ldots,t-1$, we
get, since the random sequence $\{(M_t,Z_t);t=1,2,\ldots\}$ is i.i.d.,
%
%e2.1 #&#
\begin{equation}
\label{EMarkovchainsumrep} X_t \stackrel{d} {=} \sum
_{i=1}^{t}\prod_{j=1}^{i-1}M_jZ_{i}
+ \prod_{j=1}^{t}M_{j}Z_0
\qquad\forall t=1,2,\ldots.
\end{equation}
C0 implies that $\prod_{t=1}^nM_tZ_0\xrightarrow{\mathrm{a.s.}}0$ as
$n\to\infty$. Hence, the desired conclusion follows from
(\ref{EMarkovchainsumrep}) and the Cram\'er--Slutsky theorem.

(i)${}\Rightarrow{}$(iii). C0 implies that $\prod_{t=1}^nM_tZ_0\xrightarrow{\mathrm{a.s.}}0$ as
$n\to\infty$, so by (\ref{EMarkovchainsumrep}) and the
Cram\'er--Slutsky theorem, $\sum_{i=1}^t\prod_{j=1}^{i-1}M_jZ_{i}$
converges in
distribution as $t\to\infty$. We need to prove that it also converges
a.s. We define, for brevity of notation,
\[
S_{m,n} = \sum_{i=m+1}^n\prod
_{j=1}^{i-1}M_jZ_{i}
\qquad\forall0\leq m\leq n,
\]
where $S_{n,n}=0$ for each $n\geq0$. The
following facts will be important:
%
%e2.2 #&#
\begin{equation}
\label{Eprodindepfactors} S_{m,n} = \sum_{i=m+1}^n
\prod_{j=1}^{i-1}M_jZ_{i}
= \prod_{j=1}^mM_{j}\sum
_{i=m+1}^n\prod_{j=m+1}^{i-1}M_{j}Z_{i}
\qquad \forall0\leq m< n
\end{equation}
and
%
%e2.3 #&#
\begin{equation}
\label{Esumdistridentity} \sum_{i=m+1}^n
\prod_{j=m+1}^{i-1}M_{j}Z_{i}
\stackrel{d} {=} \sum_{i=1}^{n-m}\prod
_{j=1}^{i-1}M_jZ_{i}
\qquad\forall0\leq m< n.
\end{equation}
Also, since $\sum_{i=1}^t\prod_{j=1}^{i-1}M_jZ_{i}$ converges in
distribution as $t\to\infty$, the associated sequence of
distributions is tight. Therefore, for each $\delta>0$, there exists
$K<\infty$ such that
%
%e2.4 #&#
\begin{equation}
\label{Etightnessbound} \mathbb{P}\Biggl(\Biggl|\sum_{i=1}^t
\prod_{j=1}^{i-1}M_jZ_{i}\Biggr|>K
\Biggr)<\frac{\delta}{2}\qquad\forall t=1,2,\ldots.
\end{equation}
For each
$\varepsilon>0$, each $\delta>0$, and each $n>m$, we get, if $K$ is
chosen as
in (\ref{Etightnessbound}) and $m$ is chosen
large enough,
\begin{eqnarray*}
\mathbb{P}\bigl(|S_{m,n}|>\varepsilon\bigr) &\leq& \mathbb{P}\Biggl(\Biggl\llVert
\prod_{j=1}^mM_{j}\Biggr\rrVert
>\frac
{\varepsilon}{K}\Biggr) + \mathbb{P}\Biggl(\Biggl|\sum_{i=m+1}^n
\prod_{j=m+1}^{i-1}M_{j}Z_{i}\Biggr|>K
\Biggr)
\\
&=& \mathbb{P}\Biggl(\Biggl\llVert \prod_{j=1}^mM_{j}
\Biggr\rrVert >\frac
{\varepsilon}{K}\Biggr) + \mathbb{P}\Biggl(\Biggl|\sum
_{i=1}^{n-m}\prod_{j=1}^{i-1}M_jZ_{i}\Biggr|>K
\Biggr) \leq\frac{\delta}{2} + \frac{\delta}{2}= \delta.
\end{eqnarray*}
Here, we used (\ref{Eprodindepfactors}) in the first inequality,
(\ref{Esumdistridentity}) in the equality, and C0 in the second
inequality. We conclude that
%
%e2.5 #&#
\begin{equation}
\label{Ecauchyinprobability} \sup_{n>m}\mathbb{P}\bigl(|S_{m,n}|>
\varepsilon\bigr)\to 0\mbox{ as $m\to\infty$} \qquad\forall\varepsilon>0.
\end{equation}
Our next goal is to show that, for each $\varepsilon>0$ and $m\geq0$,
if $K$ is chosen so that (\ref{Etightnessbound}) is
satisfied with $\delta=2(1-c)$, where
$0<c<1$, then:
%
%e2.6 #&#
\begin{equation}
\label{Esupremainequality} c\mathbb{P}\Bigl(\sup_{n>m}|S_{m,n}|>2
\varepsilon\Bigr) \leq\sup_{n>m}\mathbb{P}\bigl(|S_{m,n}|>
\varepsilon\bigr) + \mathbb{P}\Biggl(\bigcup_{k=m+1}^\infty
\Biggl\{\Biggl\llVert \prod_{j=1}^kM_{j}
\Biggr\rrVert >\frac{\varepsilon}{K}\Biggr\}\Biggr).
\end{equation}
To this end, we fix $\varepsilon>0$ and $m\geq0$, and note that with
this particular choice of $K$, (\ref{Esumdistridentity}) implies:
\[
\mathbb{P}\Biggl(\Biggl|\sum_{i=k+1}^n\prod
_{j=k+1}^{i-1}M_{j}Z_{i}\Biggr|
\leq K\Biggr) \geq c \qquad\forall0\leq k\leq n,
\]
which in turn gives
%
%e2.7 #&#
\begin{eqnarray}
\label{Eintermediatelowerbound}
&&\sum_{k=m+1}^{n}
\mathbb{P}\Biggl(\bigcap_{j=m+1}^{k-1}
\bigl\{|S_{m,j}|\leq 2\varepsilon\bigr\}\cap\bigl\{|S_{m,k}|>2\varepsilon\bigr\}
\Biggr)\mathbb{P}\Biggl(\Biggl|\sum_{i=k+1}^n\prod
_{j=k+1}^{i-1}M_{j}Z_{i}\Biggr|
\leq K\Biggr)
\nonumber\\[-8pt]\\[-8pt]
&&\quad\geq c\mathbb{P}\Bigl(\max_{m<k\leq
n}|S_{m,k}|>2
\varepsilon\Bigr) \qquad\forall n\geq m.
\nonumber
\end{eqnarray}
In order to obtain an upper bound for the left-hand side of
(\ref{Eintermediatelowerbound}), we note that, by the triangle inequality,
$|S_{m,k}| - |S_{k,n}| \leq|S_{m,n}|$ for each $m\leq k\leq n$. This implies:
\begin{eqnarray*}
&&\sum_{k=m+1}^{n}\mathbb{P}\Biggl(\bigcap
_{j=m+1}^{k-1}\bigl\{|S_{m,j}|\leq 2
\varepsilon\bigr\}\cap\bigl\{|S_{m,k}|>2\varepsilon\bigr\}\cap\bigl\{|S_{k,n}|\leq
\varepsilon\bigr\}\Biggr)
\\
&&\quad= \mathbb{P}\Biggl(\bigcup_{k=m+1}^{n}
\Biggl(\bigcap_{j=m+1}^{k-1}\bigl\{|S_{m,j}|
\leq 2\varepsilon\bigr\}\cap\bigl\{|S_{m,k}|>2\varepsilon\bigr\}\cap\bigl\{|S_{k,n}|
\leq \varepsilon\bigr\}\Biggr)\Biggr)
\\
&&\quad\leq\mathbb{P}\Biggl(\bigcup_{k=m+1}^n
\bigl\{|S_{m,k}|>2\varepsilon\bigr\}\cap\bigl\{ |S_{k,n}|\leq\varepsilon\bigr\}
\Biggr) \\
&&\quad\leq\mathbb{P}\bigl(|S_{m,n}|>\varepsilon\bigr)\qquad\forall n\geq m.
\end{eqnarray*}
Moreover, by (\ref{Eprodindepfactors}),
\[
\Biggl\{\Biggl\llVert \prod_{j=1}^kM_{j}
\Biggr\rrVert \leq\frac{\varepsilon
}{K}\Biggr\}\cap\Biggl\{\Biggl|\sum
_{i=k+1}^n\prod_{j=k+1}^{i-1}M_{j}Z_{i}\Biggr|
\leq K\Biggr\} \subset\bigl\{|S_{k,n}|\leq\varepsilon\bigr\}
\qquad\forall m\leq k
\leq n.
\]
Combining the last two results with the fact that the
random sequence $\{(M_t,Z_t);t=1,2,\ldots\}$ is i.i.d., we get the desired
upper bound:
\begin{eqnarray*}
&&\sum_{k=m+1}^{n}\mathbb{P}\Biggl(\bigcap
_{j=m+1}^{k-1}\bigl\{|S_{m,j}|\leq 2
\varepsilon\bigr\}\cap\bigl\{|S_{m,k}|>2\varepsilon\bigr\}\Biggr)\mathbb{P}\Biggl(\Biggl|\sum
_{i=k+1}^n\prod_{j=k+1}^{i-1}M_{j}Z_{i}\Biggr| \leq K\Biggr)
\\
&&\quad= \sum_{k=m+1}^{n}\mathbb{P}\Biggl(\bigcap
_{j=m+1}^{k-1}\bigl\{|S_{m,j}|\leq 2
\varepsilon\bigr\}\cap\bigl\{|S_{m,k}|>2\varepsilon\bigr\}\\
&&\hspace*{40pt}\qquad{}\cap\Biggl\{\Biggl|\sum
_{i=k+1}^n\prod_{j=k+1}^{i-1}M_{j}Z_{i}\Biggr|
\leq K\Biggr\}\Biggr)
\\
&&\quad= \sum_{k=m+1}^{n}\mathbb{P}\Biggl(\bigcap
_{j=m+1}^{k-1}\bigl\{|S_{m,j}|\leq 2
\varepsilon\bigr\}\cap\bigl\{|S_{m,k}|>2\varepsilon\bigr\}\cap\Biggl\{\Biggl\llVert
\prod_{j=1}^kM_{j}\Biggr\rrVert
\leq\frac{\varepsilon}{K}\Biggr\}
\\
&&\hspace*{39.7pt}\qquad{}\cap\Biggl\{\Biggl|\sum_{i=k+1}^n\prod
_{j=k+1}^{i-1}M_{j}Z_{i}\Biggr|\leq K
\Biggr\}\Biggr)
\\
&&\qquad{}+ \sum_{k=m+1}^{n}\mathbb{P}\Biggl(\bigcap
_{j=m+1}^{k-1}\bigl\{|S_{m,j}|\leq 2
\varepsilon\bigr\}\cap\bigl\{|S_{m,k}|>2\varepsilon\bigr\}\cap\Biggl\{\Biggl\llVert
\prod_{j=1}^kM_{j}\Biggr\rrVert
>\frac{\varepsilon}{K}\Biggr\}
\\
&&\hspace*{51pt}\qquad{}\cap\Biggl\{\Biggl|\sum_{i=k+1}^n\prod
_{j=k+1}^{i-1}M_{j}Z_{i}\Biggr|\leq K
\Biggr\}\Biggr)
\\
&&\quad\leq\mathbb{P}\bigl(|S_{m,n}|>\varepsilon\bigr) + \mathbb{P}\Biggl(\bigcup
_{k=m+1}^n\Biggl\{\Biggl\llVert \prod
_{j=1}^kM_{j}\Biggr\rrVert >
\frac{\varepsilon}{K}\Biggr\}\Biggr)\qquad\forall n\geq m.
\end{eqnarray*}
Letting
$n\to\infty$ (and remembering that $m\geq0$ is fixed), the last
result and (\ref{Eintermediatelowerbound}) together imply
(\ref{Esupremainequality}).

Finally, by (\ref{Esupremainequality})
and the triangle inequality,
\begin{eqnarray*}
\mathbb{P}\Bigl(\mathop{\sup_{m<k,\ell}}_{k<\ell}|S_{k,\ell
}|>4
\varepsilon\Bigr) &\leq& \mathbb{P}\Bigl(\sup_{n>m}|S_{m,n}|>2
\varepsilon\Bigr)
\\
&\leq&\frac{1}{c}\sup_{n>m}\mathbb{P}\bigl(|S_{m,n}|>
\varepsilon\bigr)\\
&&{} + \frac{1}{c}\mathbb{P}\Biggl(\bigcup
_{k=m+1}^\infty\Biggl\{\Biggl\llVert \prod
_{j=1}^kM_{j}\Biggr\rrVert >
\frac{\varepsilon}{K}\Biggr\}\Biggr) \qquad\forall\varepsilon>0, m\geq0.
\end{eqnarray*}
By (\ref{Ecauchyinprobability}), the
first term on the right-hand side converges to 0
as $m\to\infty$, while the second term converges to 0 as $m\to\infty$
by C0. Hence, $\sup_{{m<k,\ell}\atop {k<\ell}}|S_{k,\ell}|$
converges in
probability to 0 as
\mbox{$m\to\infty$}. However, by definition,
$\sup_{{m<k,\ell}\atop {k<\ell}}|S_{k,\ell}|$ decreases
monotonically a.s. to a nonnegative random variable as
$m\to\infty$. To avoid a contradiction, this random variable must be
0 with
probability 1. It follows that, with
probability 1,
$\{\sum_{i=1}^t\prod_{j=1}^{i-1}M_jZ_{i};t=1,2,\ldots\}$ is a
Cauchy sequence, so $\lim_{t\to\infty}\sum_{i=1}^t\prod_{j=1}^{i-1}M_jZ_{i}$
exists a.s.

(ii)${}\Rightarrow{}$(iii)${}\Rightarrow{}$(iv)${}\Rightarrow{}$(v). Immediate.

(iv)${}\Rightarrow{}$(vi). As stated in
Remark \ref{Rimplicationswithoutnormconv}, C0 is not needed
to prove this
implication. Instead, we use the theorem in \cite{KS64}, also known as the
\emph{Kochen--Stone
lemma}. By this theorem (or lemma), for any sequence of
events $\{A_t;t=1,2,\ldots\}$ such that $\sum_{t=1}^\infty\mathbb
{P}(A_t)=\infty$ and\vspace*{-1pt}
%
%e2.8 #&#
\begin{equation}
\label{Ekochenstonecondition} \limsup_{n\to\infty}\frac{ (\sum_{t=1}^n\mathbb{P}(A_t)
)^2}{\sum_{r=1}^n\sum_{t=1}^n\mathbb{P}(A_r\cap
A_t)} =
c >0,
\end{equation}
it holds that $\mathbb{P}(A_t\mbox{ i.o.})\geq c$. Define the random
sequence $\{Y_t;t=1,2,\ldots\}$ by:\vspace*{-1pt}
\[
Y_t=\min_{k=1,\ldots,t-1}\Biggl|\prod_{j=k}^{t-1}M_jZ_t\Biggr|
\qquad\forall t=1,2,\ldots.
\]
Recall that by definition $Y_1=\infty$ (since it is the
minimum over an empty set). Let $x>0$, and define the events
$\{A_t;t=1,2,\ldots\}$ by: $A_t=\{Y_t>x\}$ $\forall t=1,2,\ldots\,$. We
note that if (vi) does not hold, then
$\sum_{t=1}^\infty\mathbb{P}(A_t)=\infty$ for some $x>0$. We will show
that in this case (\ref{Ekochenstonecondition}) holds with
$c\geq\frac{1}{2}$, implying that
$\mathbb{P}(|\prod_{j=1}^{t-1}M_jZ_t|>x \mbox{ i.o.}) \geq
\mathbb{P}(Y_t>x \mbox{ i.o.})\geq\frac{1}{2}>0$. Hence, (iv) does not
hold.

For the probabilities in
the denominator of (\ref{Ekochenstonecondition}), we get, if $1\leq r<t$,\vspace*{-1pt}
\begin{eqnarray*}
\mathbb{P}\bigl(\{Y_r> x\}\cap\{Y_t> x\}\bigr) &=&
\mathbb{P}\Biggl(\{Y_r> x\}\cap\Biggl\{ \min_{k=1,\ldots,t-1}\Biggl|
\prod_{j=k}^{t-1}M_jZ_t\Biggr|>
x\Biggr\}\Biggr)
\\
&\leq& \mathbb{P}(Y_r> x)\mathbb{P}\Biggl(\min_{k=r+1,\ldots,t-1}\Biggl|
\prod_{j=k}^{t-1}M_jZ_t\Biggr|>
x\Biggr) \\
&=& \mathbb{P}(Y_r> x)\mathbb{P}(Y_{t-r}> x).
\end{eqnarray*}
This implies that
\begin{eqnarray*}
&&
\sum_{r=1}^n\sum
_{t=1}^n\mathbb{P}\bigl(\{Y_r> x\}\cap
\{Y_t> x\}\bigr)
\\
&&\quad\leq\sum_{r=1}^n\mathbb{P}(Y_r>
x) + 2\sum_{r=1}^{n-1}\mathbb
{P}(Y_r> x)\sum_{t=r+1}^n
\mathbb{P}(Y_{t-r}> x)
\\
&&\quad\leq\sum_{r=1}^n\mathbb{P}(Y_r>
x) + 2\sum_{r=1}^{n}\mathbb
{P}(Y_r> x)\sum_{s=1}^{n}
\mathbb{P}(Y_s> x).
\end{eqnarray*}
Hence, we obtain:
\begin{eqnarray*}
&&
\limsup_{n\to\infty}\frac{ (\sum_{t=1}^n\mathbb{P}(Y_t>
x) )^2}{\sum_{r=1}^n\sum_{t=1}^n\mathbb{P}(\{Y_r> x\}\cap\{
Y_t> x\})}
\\
&&\quad\geq\lim_{n\to\infty}\frac{ (\sum_{t=1}^n\mathbb{P}(Y_t>
x) )^2}{\sum_{r=1}^n\mathbb{P}(Y_r> x)
+ 2 (\sum_{t=1}^n\mathbb{P}(Y_t> x) )^2} = \frac{1}{2}.
\end{eqnarray*}

(vi)${}\Rightarrow{}$(ii). This part of the proof is divided into several
steps. First, we prove that if $\llVert {\prod_{t=1}^nM_t}\rrVert \xrightarrow{\mathrm{a.s.}}0$ as
$n\to\infty$, then
%
%e2.9 #&#
\begin{equation}
\label{Enormofproductscondition} \sum_{t=1}^\infty
\mathbb{P}\Biggl(\min_{k=1,\ldots,t-1}\Biggl\llVert \prod
_{j=k}^{t-1}M_j\Biggr\rrVert > x\Biggr)<
\infty \qquad\forall x>0.
\end{equation}
We use the Kochen--Stone lemma, as in the preceding part of the proof.
Let
\[
U_t=\min_{k=1,\ldots,t-1}\Biggl\llVert {\prod_{j=k}^{t-1}M_j}\Biggr\rrVert
\qquad\forall t=1,2,\ldots.
\]
Let $x>0$, and define the events $\{
A_t;t=1,2,\ldots\}$ by: $A_t=\{U_t>x\}$ $\forall t=1,2,\ldots\,$. Assume
that $\sum_{t=1}^\infty\mathbb {P}(A_t)=\infty$. As before, for the
probabilities in the denominator of (\ref{Ekochenstonecondition}), we
get:
\[
\mathbb{P}\bigl(\{U_r> x\}\cap\{U_t> x\}\bigr) \leq
\mathbb{P}(U_r> x)\mathbb{P}(U_{t-r}> x) \qquad\forall1\leq r<t,
\]
implying that
\[
\limsup_{n\to\infty}\frac{ (\sum_{t=1}^n\mathbb{P}(U_t>
x) )^2}{\sum_{r=1}^n\sum_{t=1}^n\mathbb{P}(\{U_r>
x\}\cap\{U_s> x\})} \geq\frac{1}{2},
\]
so $\mathbb{P}(U_t> x \mbox{ i.o.})\geq\frac{1}{2}$. Hence, it cannot
hold that $\llVert {\prod_{t=1}^nM_t}\rrVert
\xrightarrow{\mathrm{a.s.}}0$ as $n\to\infty$.

Next, let as before $Y_t=\min_{k=1,\ldots,t-1}|{\prod_{j=k}^{t-1}M_jZ_t}|$
$\forall t=1,2,\ldots\,$. Since
\[
\Biggl|\prod_{j=1}^{t-1}M_jZ_t\Biggr|
\leq \Biggl\llVert \prod_{j=1}^{k-1}M_{j}
\Biggr\rrVert \Biggl|\prod_{j=k}^{t-1}M_jZ_t\Biggr|
\qquad\forall t=1,2,\ldots; k=1,\ldots,t-1,
\]
it holds that
\[
\Biggl|\prod_{j=1}^{t-1}M_jZ_t\Biggr|
\leq\sup_{n\geq
0}\Biggl\llVert \prod_{i=1}^nM_i
\Biggr\rrVert \min_{k=1,\ldots,t-1}\Biggl|\prod_{j=k}^{t-1}
M_jZ_t\Biggr| \qquad\forall t=1,2,\ldots,
\]
where, since\vspace*{2pt} $\llVert {\prod_{t=1}^nM_t}\rrVert \xrightarrow
{\mathrm{a.s.}}0$ as
$n\to\infty$, $\sup_{n\geq0}\llVert {\prod_{i=1}^nM_i}\rrVert <\infty$
a.s. This implies that in order to prove (ii), it is sufficient to
prove that $\sum_{t=1}^\infty Y_t <\infty$ a.s.

Furthermore, by Fubini's theorem,
%
%e2.10 #&#
\begin{eqnarray}
\label{EKolmogorovexpectationrepresentation} \mathbb{E}\bigl(Y_tI
\{Y_t\leq1\}\bigr) &=& \int_{(0,1]}y\mrmd F_{Y_t}(y)
\nonumber\\
&=& \int_0^{1}\mathbb{P}(x<Y_t
\leq1)\mrmd x \\
&\leq&\int_0^1\mathbb{P}(Y_t>x)\mrmd x
\qquad\forall t=1,2,\ldots,
\nonumber
\end{eqnarray}
implying that
%
%e2.11 #&#
\begin{equation}
\label{EKolmogorovothercondition} \sum_{t=1}^\infty
\mathbb{E}\bigl(Y_tI\{Y_t\leq1\}\bigr) \leq \int
_0^1\sum_{t=1}^\infty
\mathbb{P}(Y_t>x)\mrmd x.
\end{equation}
We note that, by (vi), $\sum_{t=1}^\infty\mathbb{P}(Y_t>x)<\infty$ for
each $x>0$. We will prove that the right-hand side of
(\ref{EKolmogorovothercondition}) is finite. By monotone convergence,
this will imply that
\[
\mathbb{E}\Biggl(\sum_{t=1}^\infty
Y_tI\{Y_t\leq1\}\Biggr) = \sum
_{t=1}^\infty\mathbb{E}\bigl(Y_tI
\{Y_t\leq1\}\bigr) < \infty,
\]
from which it will follow that $\sum_{t=1}^\infty Y_tI\{Y_t\leq
1\}<\infty$ a.s. Since, by (vi) and the Borel--Cantelli lemma,
$Y_t\xrightarrow{\mathrm{a.s.}}0$ as
$t\to\infty$, we will be able to conclude that $\sum_{t=1}^\infty
Y_t<\infty$ a.s.

Define $\{\widetilde Y_t;t=1,2,\ldots\}$ by $\widetilde Y_t=
\min_{k=1,\ldots,t-1}|\prod_{j=k}^{t-1}M_{t-j}\widetilde Z_1|$ $\forall
t=1,2,\ldots\,$, where $\widetilde Z_1$ is a random variable independent
of $\{(M_t,Z_t);t=1,2,\ldots\}$ such that $\widetilde Z_1
\stackrel{d}{=} Z_1$. By definition, $\{\widetilde Y_t;t=1,2,\ldots\}$
is a nonincreasing random sequence, while clearly also $\widetilde Y_t
\stackrel{d}{=} Y_t$ $\forall t=1,2,\ldots$ (in particular, $\widetilde
Y_1 = Y_1 = \infty$, since they are both minima over empty sets).
Define, for each $x>0$, the random variable
\[
T_x=\inf\{t=1,2,\ldots;\widetilde Y_t\leq x\} = \inf
\Biggl\{t=1,2,\ldots ;\Biggl|\prod_{j=1}^{t-1}M_{t-j}
\widetilde Z_1\Biggr|\leq x\Biggr\}.
\]
Clearly, $T_x$ is a stopping time with respect to the filtration
$\{\mathscr{G}_t;t=1,2,\ldots\}$, defined by: $\mathscr{G}_t =
\sigma(\widetilde Z_1;M_1,\ldots,M_{t-1})$ $\forall t=1,2,\ldots\,$. Moreover,
%
%e2.12 #&#
\begin{equation}
\label{Eexpectedhittingtime} \sum_{t=1}^\infty
\mathbb{P}(Y_t>x) = \sum_{t=1}^\infty
\mathbb{P}(\widetilde Y_t>x) = \sum_{t=1}^\infty
\mathbb{P}(T_x>t) = \mathbb{E}(T_x)-1,
\end{equation}
so (vi) implies that $\mathbb{E}(T_x)<\infty$ for each $x>0$. Define,
for each $x>0$, the random variables $T_x^{(1)}=T_1$ and
\[
T_x^{(2)}=\inf\Biggl\{t=1,2,\ldots;\Biggl\llVert \prod
_{j=1}^{t}M_{T_1+t-j}\Biggr\rrVert
\leq x\Biggr\}.
\]
Since
$\{M_t;t=1,2,\ldots\}$ are i.i.d. and independent of $\widetilde Z_1$,
it holds
that $\{M_s;s=t,t+1,\ldots\}$ are independent of $\mathscr{G}_t$ for
each $t=1,2,\ldots\,$. Since $T_1$ is an a.s. finite stopping time with
respect to
$\{\mathscr{G}_t;t=1,2,\ldots\}$, we get:
\begin{eqnarray*}
&&
\mathbb{P}\bigl(\bigl\{T_x^{(2)}>t\bigr\}\cap
\{T_1=r\}\bigr)\\
&&\quad= \mathbb{P}\Biggl(\Biggl\{\min_{k=1,\ldots,t}
\Biggl\llVert \prod_{j=k}^{t}M_{T_1+t-j}
\Biggr\rrVert >x\Biggr\}\cap\{T_1=r\}\Biggr)
\\
&&\quad= \mathbb{P}\Biggl(\Biggl\{\min_{k=1,\ldots,t}\Biggl\llVert \prod
_{j=k}^{t}M_{r+t-j}\Biggr\rrVert >x
\Biggr\}\cap\{T_1=r\}\Biggr)
\\
&&\quad= \mathbb{P}\Biggl(\min_{k=1,\ldots,t}\Biggl\llVert \prod
_{j=k}^{t}M_{j}\Biggr\rrVert >x\Biggr)
\mathbb{P}(T_1=r) \qquad\forall t=1,2,\ldots; r=1,2,\ldots.
\end{eqnarray*}
In particular,
\[
\mathbb{P}\bigl(T_x^{(2)}>t\bigr) = \mathbb{P}\Biggl(\min
_{k=1,\ldots,t}\Biggl\llVert \prod_{j=k}^{t}M_{j}
\Biggr\rrVert >x\Biggr) \qquad\forall t=1,2,\ldots
\]
and
\[
\mathbb{E}\bigl(T_x^{(2)}\bigr)-1 = \sum
_{t=1}^\infty\mathbb{P}\bigl(T_x^{(2)}>t
\bigr) = \sum_{t=1}^\infty\mathbb{P}\Biggl(\min
_{k=1,\ldots,t}\Biggl\llVert \prod_{j=k}^{t}M_{j}
\Biggr\rrVert >x\Biggr) < \infty\qquad\forall x>0,
\]
where finiteness follows from (\ref{Enormofproductscondition}).

Repeating this process, we define recursively, for each $x>0$, the
random variables
$\{T_x^{(k)};k=2,3,\ldots\}$ by:
\[
T_x^{(k)}=\inf\Biggl\{t=1,2,\ldots;\Biggl\llVert \prod
_{j=1}^{t}M_{S_x^{(k-1)}+t-j}\Biggr\rrVert
\leq x\Biggr\} \qquad\forall k=2,3,\ldots,
\]
where $S_x^{(k)}=\sum_{i=1}^kT_x^{(i)}$ $\forall
k=1,2,\ldots\,$. Since $\{M_s;s=t,t+1,\ldots\}$ are independent of
$\mathscr{G}_t$ for
each $t=1,2,\ldots\,$, and since $\{S_x^{(k)};k=1,2,\ldots\}$ are
stopping times with respect to
$\{\mathscr{G}_t;t=1,2,\ldots\}$, we see that
$\{T_x^{(k)};k=2,3,\ldots\}$ are i.i.d. with finite mean.

We now observe that by the submultiplicative property,
\[
\Biggl|\prod_{j=1}^{S_x^{(k+1)}-1}M_{S_x^{(k+1)}-j}
\widetilde Z_1\Biggr| \leq \Biggl|\prod_{j=1}^{T_1-1}M_{T_1-j}
\widetilde Z_1\Biggr|\prod_{i=2}^{k+1}
\Biggl\llVert \prod_{j=1}^{T_x^{(i)}}M_{S_x^{(i)}-j}
\Biggr\rrVert \leq x^k \qquad\forall k=1,2,\ldots; x>0,
\]
which implies that
\[
T_x \leq S_{x^{1/k}}^{(k+1)} = T_1 +
T_{x^{1/k}}^{(2)} + \cdots+ T_{x^{1/k}}^{(k+1)}
\qquad\forall k=1,2,\ldots; x>0.
\]
Taking expectations on both sides in this inequality gives:
\[
\mathbb{E}(T_x) \leq\mathbb{E}(T_1) + k\mathbb{E}
\bigl(T_{x^{1/k}}^{(2)}\bigr) \qquad\forall k=1,2,\ldots; x>0.
\]
Choosing $a\in(0,1)$ and letting
$k_x=\lceil\frac{\log x}{\log a}\rceil$ $\forall x\in(0,1)$, we get:
\[
x^{1/k_x} = \exp \biggl(\frac{\log x}{\lceil{\log x}/{\log
a}\rceil} \biggr) \geq a \qquad\forall x
\in(0,1),
\]
implying that
\[
\mathbb{E}(T_x) \leq\mathbb{E}(T_1) + k_x
\mathbb{E}\bigl(T_a^{(2)}\bigr) \leq\mathbb{E}(T_1)
+ \mathbb{E}\bigl(T_a^{(2)}\bigr) \biggl(\frac{\log x}{\log a}
+ 1\biggr) \qquad\forall x\in(0,1).
\]
This combined with (\ref{Eexpectedhittingtime}) implies that the
right-hand side of (\ref{EKolmogorovothercondition}) is finite, since
\[
\int_0^1\log x\mrmd x = \lim_{\epsilon\to0}
\int_\epsilon^1\log x\mrmd x = \lim_{\epsilon\to0}[x
\log x -x]_\epsilon^1 = -1.
\]

(v)${}\Rightarrow{}$(iv). Since
\[
\Biggl|\prod_{j=1}^{t-1}M_jZ_t\Biggr|
\leq \Biggl\llVert \prod_{j=1}^{m-1}M_j
\Biggr\rrVert \Biggl|\prod_{j=m}^{t-1}M_jZ_t\Biggr|
\qquad\forall1\leq m\leq t,
\]
it holds for each $\varepsilon>0$ and $K>0$ that
\begin{eqnarray*}
\mathbb{P}\Biggl(\bigcup_{t=m}^n\Biggl
\{\Biggl|\prod_{j=1}^{t-1}M_jZ_t\Biggr|>
\varepsilon\Biggr\}\Biggr) &\leq&\mathbb{P}\Biggl(\Biggl\llVert \prod
_{j=1}^{m-1}M_j\Biggr\rrVert >
\frac{\varepsilon}{K}\Biggr)
\\
&&{}+ \mathbb{P}\Biggl(\bigcup_{t=m}^n\Biggl
\{\Biggl|\prod_{j=m}^{t-1}M_jZ_t\Biggr|>K
\Biggr\}\Biggr) \qquad \forall1\leq m\leq n.
\end{eqnarray*}
For the second term on the right-hand side, since the random sequence
$\{(M_t,Z_t);t=1,2,\ldots\}$ is i.i.d.,
\begin{eqnarray*}
\mathbb{P}\Biggl(\bigcup_{t=m}^n\Biggl
\{\Biggl|\prod_{j=m}^{t-1}M_jZ_t\Biggr|>K
\Biggr\}\Biggr) &=& \mathbb{P}\Biggl(\bigcup_{t=m}^n
\Biggl\{\Biggl|\prod_{j=m}^{t-1}M_{j-m+1}Z_{t-m+1}\Biggr|>K
\Biggr\}\Biggr)
\\
&=& \mathbb{P}\Biggl(\bigcup_{t=m}^n\Biggl
\{\Biggl|\prod_{j=1}^{t-m}M_jZ_{t-m+1}\Biggr|>K
\Biggr\}\Biggr) \\
&=& \mathbb{P}\Biggl(\bigcup_{t=1}^{n-m+1}
\Biggl\{\Biggl|\prod_{j=1}^{t-1}M_jZ_t\Biggr|>K
\Biggr\}\Biggr)
\\
&\leq&\mathbb{P}\Biggl(\sup_{t=1,2,\ldots}\Biggl|\prod
_{j=1}^{t-1}M_jZ_t\Biggr|>K\Biggr)
\qquad\forall1\leq m\leq n.
\end{eqnarray*}
Fixing $m\geq1$ and letting $n\to\infty$, we get:
\begin{eqnarray*}
\mathbb{P}\Biggl(\bigcup_{t=m}^\infty\Biggl
\{\Biggl|\prod_{j=1}^{t-1}M_jZ_t\Biggr|>
\varepsilon \Biggr\}\Biggr) &\leq&\mathbb{P}\Biggl(\Biggl\llVert \prod
_{j=1}^{m-1}M_j\Biggr\rrVert >
\frac{\varepsilon}{K}\Biggr)
\\
&&{}+ \mathbb{P}\Biggl(\sup_{t=1,2,\ldots}\Biggl|\prod
_{j=1}^{t-1}M_jZ_t\Biggr|>K\Biggr)
\qquad\forall m\geq1.
\end{eqnarray*}
For each $\delta>0$, by (v), the second term on the right-hand side
can be made less than $\frac{\delta}{2}$ by choosing $K$ large
enough. Similarly, using C0, the first term on the right-hand side
can be made less than $\frac{\delta}{2}$ by choosing $m$ large
enough. This gives:
\[
\mathbb{P}\Biggl(\bigcup_{t=m}^\infty\Biggl
\{\Biggl|\prod_{j=1}^{t-1}M_jZ_t\Biggr|>
\varepsilon \Biggr\}\Biggr) \leq\frac{\delta}{2} + \frac{\delta}{2} = \delta,
\]
which implies (iv).
\end{pf*}

%s3 #&#
\section{Counterexamples and special cases} \label{Scounterexamples}
In this section, we consider some counterexamples, some
special cases, and a condition on the matrices $\{M_t;t=1,2,\ldots\}$
which is only sufficient for C0, but
somewhat easier to validate. In
Example \ref{RC0necessarycounterexample}, it is shown that in
the case $d>1$, (ii) in
Theorem \ref{Tmultidimperpetuities} does not imply C0. In
Examples \ref{Riitoicounterexample}--\ref{Rvtovicounterexample}, it is shown that in the case $d>1$, if C0 does not
hold, not all of the conclusions of
Theorem \ref{Tmultidimperpetuities} hold. The
special cases considered are the case $d=1$ (completely solved in \cite
{GM00}), and the case when
$M_t=M$ $\forall t=1,2,\ldots\,$, where $M$ is a (deterministic)
constant matrix.

%ex3.1 #&#
\begin{example} \label{RC0necessarycounterexample}
Consider first the case $d=1$. This case was completely solved in
\cite{GM00}, where it was shown that if $\mathbb{P}(Z_1=0)<1$, then
(vi) implies C0, and if
also $\mathbb{P}(|M_1|=1)<1$, then (v) implies C0. Moreover, if
$\mathbb{P}(Z_1=0)<1$, then clearly (iv) implies that
$\mathbb{P}(|M_1|=1)<1$. As a consequence, if
$d=1$ and $\mathbb{P}(Z_1=0)<1$, then (ii), (iii), (iv), (v)
\emph{combined with $\mathbb{P}(|M_1|=1)<1$}, and (vi) are
equivalent, and they
all imply C0.

However, if $d>1$, the following counterexample shows
that even if $\mathbb{P}(Z_1=0)<1$, (ii) does not imply C0. Let $d=2$,
and let $v_1$
and $v_2$ be orthonormal column vectors in $\mathbb{R}^2$. Let
$0<\alpha<1$. Define $M_t=\alpha v_1v_1^T + v_2v_2^T$ $\forall
t=1,2,\ldots\,$, and $Z_t=v_1$ $\forall
t=1,2,\ldots\,$. Then, $\prod_{j=1}^{t-1}M_jZ_t = \alpha^{t-1}v_1$
$\forall t=1,2,\ldots\,$, so (ii) holds. On the other hand,
$\llVert {\prod_{j=1}^{t}M_j}\rrVert = 1$
$\forall t=1,2,\ldots\,$, which does not
converge to 0 a.s. as $t\to\infty$.
\end{example}

%ex3.2 #&#
\begin{example} \label{Riitoicounterexample}
If $d>1$ and C0 does not hold,
then the implication (ii)${}\Rightarrow{}$(i) does not hold. To see this,
let $d=2$, and let $v_1$
and $v_2$ be orthonormal column vectors in $\mathbb{R}^2$. Let
$0<\alpha<1<\beta<\infty$. Define $M_t=\alpha v_1v_1^T + \beta
v_2v_2^T$ $\forall t=1,2,\ldots\,$, and $Z_t=v_1$ $\forall
t=1,2,\ldots\,$. Let $Z_0=v_2$. Then, $\prod_{j=1}^{t-1}M_jZ_t = \alpha
^{t-1}v_1$ $\forall t=1,2,\ldots\,$, so (ii) holds, and
$\sum_{i=1}^t\prod_{j=1}^{i-1}M_jZ_t$ converges a.s. to
$\frac{1}{1-\alpha}v_1$ (a\vspace*{1pt} deterministic vector) as $t\to\infty$.
On the other
hand, $\llVert {\prod_{j=1}^{t}M_j}\rrVert = \beta^t$
$\forall t=1,2,\ldots\,$, which does not
converge to 0 a.s. as $t\to\infty$. If (i) holds, then by
(\ref{EMarkovchainsumrep}), (ii) and the Cram\'er--Slutsky theorem,
$\prod_{j=1}^{t}M_jZ_0$
must converge in distribution as $t\to\infty$. However, $\prod_{j=1}^{t}M_jZ_0 = \beta^tv_2$
$\forall t=1,2,\ldots\,$, which does not converge in distribution as
$t\to\infty$ (the corresponding sequence of distributions is not
tight). Hence, (i) does not hold.
\end{example}

%ex3.3 #&#
\begin{example} \label{Ritovandvicounterexample}
If C0 does not hold, then the implications (i)${}\Rightarrow{}$(v) and
(i)${}\Rightarrow{}$(vi) do not hold. To see this, let $d=1$, $|\beta|>1$ and
$c>0$. Define $M_t=\beta$ $\forall t=1,2,\ldots\,$, $Z_t=(1-\beta)c$
$\forall
t=1,2,\ldots\,$, and $Z_0=c$. (This is an example where the
``nondegeneracy'' condition
(2.7) in \cite{GM00} does not hold.) Then
\[
\sum_{i=1}^t\prod
_{j=1}^{i-1}M_jZ_t + \prod
_{j=1}^{t}M_jZ_0 =
(1-\beta)c\frac{1-\beta^t}{1-\beta} + c\beta^t = c \qquad\forall t=1,2,\ldots,
\]
so by (\ref{EMarkovchainsumrep}) (i) holds. On the other
hand, $\llVert {\prod_{j=1}^{t}M_j}\rrVert = |\beta|^t$
$\forall t=1,2,\ldots\,$, which does not
converge to 0 a.s. as $t\to\infty$. Also, $|{\prod_{j=1}^{t-1}M_jZ_t}|
=|(1-\beta)|c|\beta|^{t-1}$
$\forall t=1,2,\ldots\,$, so neither (v) nor (vi) holds.
\end{example}

%ex3.4 #&#
\begin{example} \label{Rvtovicounterexample}
If $d>1$ and C0 does
not hold, then the implication (v)${}\Rightarrow{}$(vi) does not
hold. To see this, we use the same setup as in
Example \ref{RC0necessarycounterexample}, except that we now define
$Z_t=v_2$ $\forall
t=1,2,\ldots\,$. Then, $\prod_{j=1}^{t-1}M_jZ_t = v_2$
$\forall t=1,2,\ldots\,$, so (v) holds, but not (vi). Moreover,
$\llVert {\prod_{j=1}^{t}M_j}\rrVert = 1$
$\forall t=1,2,\ldots\,$.
\end{example}

%re3.1 #&#
\begin{remark}[(An open problem)] \label{Rvitovcounterexample}
Despite some effort, we have not been able to find a
counterexample showing that if $d>1$
and C0 does not hold, the implication (vi)${}\Rightarrow{}$(v) does not
hold. It is therefore possible that, if $d>1$, even when C0 does not
hold, (vi) implies one or
several of (ii), (iii), (iv) or (v). We leave it as an open problem to
prove these assertions, or to disprove them by means of counterexamples.
\end{remark}

%re3.2 #&#
\begin{remark} \label{Rproductofnormscounterexample}
Consider again the case $d>1$. As pointed out in Remark 2.13 in
\cite{GM00}, a \emph{sufficient} condition for (ii) to hold is that
$\sum_{t=1}^\infty\prod_{j=1}^{t-1}\llVert  M_j\rrVert
|Z_{t}|<\infty$
a.s. By Theorem 2.1 in
\cite{GM00} (see also Example \ref{RC0necessarycounterexample}
above), the latter condition is equivalent to
\[
\sum_{t=1}^\infty\mathbb{P}\Biggl(\min
_{k=1,\ldots,t-1}\prod_{j=k}^{t-1}
\llVert M_j\rrVert |Z_t|> x\Biggr)<\infty \qquad\forall
x>0
\]
and to $\prod_{j=1}^{t-1}\llVert
M_j\rrVert |Z_{t}|\xrightarrow{\mathrm{a.s.}}0$ as $t\to\infty$. If
$\mathbb{P}(Z_1=0)<1$, these
equivalent conditions all imply that $\prod_{j=1}^t\llVert
M_j\rrVert \xrightarrow{\mathrm{a.s.}}0$ as
$n\to\infty$, which clearly implies C0.

However, C0 does not imply that $\prod_{j=1}^t\llVert
M_j\rrVert \xrightarrow{\mathrm{a.s.}}0$ as $t\to\infty$, as the
following counterexample shows. Let
$d=2$, and let $v_1$ and $v_2$ be orthonormal column vectors in
$\mathbb{R}^2$. Let $\{\alpha_t;t=1,2,\ldots\}$
be an i.i.d. random sequence such that $\mathbb{P}(\alpha_t=1) =
\mathbb{P}(\alpha_t=\frac{1}{2}) =
\frac{1}{2}$ $\forall t=1,2,\ldots\,$, and let
$K_t=\sum_{j=1}^tI\{\alpha_j=1\}$ $\forall t=1,2,\ldots\,$. Define
$M_t=\alpha_tv_1v_1^T + (\frac{3}{2}-\alpha_t)v_2v_2^T$ $\forall
t=1,2,\ldots\,$. Then
\[
\Biggl\llVert \prod_{j=1}^{t}M_j
\Biggr\rrVert = \max \biggl(\frac{1}{2^{t-K_t}},\frac{1}{2^{K_t}} \biggr) \qquad
\forall t=1,2,\ldots.
\]
By the second Borel--Cantelli lemma,
$K_t\xrightarrow{\mathrm{a.s.}}\infty$ and
$t-K_t\xrightarrow{\mathrm{a.s.}}\infty$ as $t\to\infty$, implying that
$\llVert {\prod_{j=1}^tM_j}\rrVert \xrightarrow{\mathrm{a.s.}}0$ as
$t\to\infty$. On the other hand, $\prod_{j=1}^t\llVert  M_j\rrVert = 1$
$\forall t=1,2,\ldots\,$.
\end{remark}

%re3.3 #&#
\begin{remark} \label{RC0sufficientconditions}
As noted in Remark \ref{Rproductofnormscounterexample}, the condition
$\prod_{j=1}^t\llVert  M_j\rrVert \xrightarrow{\mathrm{a.s.}}0$ as
$t\to\infty$ implies C0. By Proposition 2.6 in \cite{GM00} (see also
Section 4 in \cite{GM00}), the former condition
holds \emph{if and only if} one of the following two conditions hold:
\begin{eqnarray*}
&&\mbox{\textup{\hphantom{i}(i)}}\quad \mathbb{E}\bigl(\bigl|\log\llVert M_1\rrVert \bigr|\bigr)<
\infty\quad\mbox{and}\quad \mathbb{E}\bigl(\log\llVert M_1\rrVert \bigr) < 0;
\\
&&\mbox{(ii)}\quad \mathbb{E}\bigl(\log^-\llVert M_1\rrVert \bigr)=\infty
\quad\mbox{and}\quad \mathbb{E} \biggl(\frac{\log^+\llVert  M_1\rrVert }{A_M(\log^+\llVert
M_1\rrVert )} \biggr)<\infty,
\end{eqnarray*}
where $A_M(y)=\int_0^y\mathbb{P}(-\log\llVert  M_1\rrVert
>x)\mrmd x$ $\forall
y>0$, $\log^+x = \log(x\vee1)$ $\forall x>0$, and $\log^-x = -\log
(x\wedge1)$ $\forall x>0$.
\end{remark}

%re3.4 #&#
\begin{remark} \label{RLyapunovexponents}
Under the condition $\mathbb{E}({\log^+}\llVert  M_1\rrVert
)<\infty$,
Kingman's subadditive ergodic theorem can be used to show that
\[
\frac{1}{t}\log\Biggl\llVert \prod_{j=1}^tM_j
\Biggr\rrVert \xrightarrow {\mathrm{a.s.}}\lambda= \lim_{n\to\infty}
\frac{1}{n}\mathbb{E}\Biggl(\log\Biggl\llVert \prod
_{j=1}^nM_j\Biggr\rrVert \Biggr)\qquad
\mbox{as $t\to\infty$},
\]
where $\lambda\in[-\infty,\infty)$ is a deterministic
constant; see Theorem 6 in \cite{Kin73} and Theorem 2 in
\cite{FK60}. (Recall that the matrix norm used in these papers is
equivalent to the spectral norm.) The constant $\lambda$ is sometimes
called the
\emph{maximal Lyapunov exponent}. In particular, if $\mathbb{E}(\log
^+\llVert
M_1\rrVert )<\infty$, then C0 holds if $\lambda<0$, and does not hold
if $\lambda>0$. For more information, see \cite{FK60,Kin73}
and the references therein.
\end{remark}

%re3.5 #&#
\begin{remark} \label{Rconstantcoefficient}
Finally, consider the case when $\mathscr{L}(M_1)$ is degenerate at a constant
$d\times d$-matrix $M$, that is, the case when the RCA(1) process
$\{X_t;t=1,2,\ldots\}$ is an AR(1)
process. In this case, $\prod_{j=1}^tM_j=M^t$ $\forall t=1,2,\ldots\,$,
and the following \emph{spectral representation} holds:
%
%e3.1 #&#
\begin{equation}
\label{Espectralrepresentation} M^t = \sum
_{k=1}^s\sum_{j=0}^{m_k-1}
\biggl[\frac{d^j}{dx^j}x^t \biggr]_{x=\lambda_k}Z_{k,j}
\qquad\forall t=1,2,\ldots,
\end{equation}
where $\{\lambda_k;k=1,\ldots,s\}$ are the distinct eigenvalues of
$M$, and $\{m_k;k=1,\ldots,s\}$ are the multiplicities (all
positive integers) of the eigenvalues as zeros of the \emph{minimal
annihilating polynomial} of
$M$. Moreover, $\{Z_{k,j};k=1,\ldots,s;j=0,\ldots,m_k-1\}$ are linearly
independent $d\times d$-matrices called the \emph{components} of $M$;
for more information, see Section 9.5 in \cite{LT85}. Assuming that
$\lambda_1$ is an eigenvalue of maximum modulus, there are
two possible cases. If $|\lambda_1|<1$,
then, applying the triangle inequality to the right-hand side of
(\ref{Espectralrepresentation}), we see that $\llVert  M^t\rrVert \to0$ as
$t\to\infty$. On the other hand, if $|\lambda_1|\geq1$, then
$\llVert  M^t\rrVert \geq|M^tv_1| = |\lambda_1|^t\geq1$
$\forall
t=1,2,\ldots\,$, where $v_1$ is a normalized eigenvector corresponding
to $\lambda_1$. Hence, C0 holds if and only if $|\lambda_1|<1$.
\end{remark}

%s4 #&#
\section{Suggestions for future research} \label{Sfurtherwork}
We mention two possible research directions. First, the open problem
stated in
Remark \ref{Rvitovcounterexample}: to determine whether, in the case
$d>1$, (vi) in Theorem \ref{Tmultidimperpetuities} implies one or
several of (ii), (iii), (iv) or (v), without condition C0 (or
replacing C0 with an even less restrictive condition).
Second, to find a natural generalization (if it exists) of
the integral condition (2.1) in Theorem 2.1 in \cite{GM00} to higher
dimensions.

% zodis "Acknowledgments" paliekamas pagal autoriu

%suskaldyti doi

% imsref loaded by lrinkeviciute, 2013-06-04 08:40:13

\printhistory

\end{document}